\documentclass[a4paper,12pt]{article}

\usepackage[top=1in, bottom=1in, left=1in, right=1in]{geometry}
\usepackage{amsmath, amssymb}
\usepackage{amsthm}
\usepackage{dsfont}
\usepackage{enumerate}
\usepackage{tikz-cd}
\usepackage{graphicx}
\usepackage{spline}

\title{An upper bound on the regularity of the first homology of spline complexes}
\author{Beihui Yuan}
\begin{document}
\maketitle
\begin{abstract}
Let $\Delta$ be a connected, pure $2$-dimensional simplicial complex embedded in $\RR^2$ and let $C^{r}(\hat{\Delta})$ be the homogenized spline module of $\Delta$ with smoothness $r$ as in \cite{schenck1997local}. To study $C^{r}(\hat{\Delta})$, Schenck and Stillman developed in \cite{schenck1997local} the spline complex $S_\bullet/J_\bullet$. In \cite{schenck2002cohomology}, Schenck and Stiller conjectured that the regularity of $H_1(S_\bullet/J_\bullet)$ is less than $2r+1$. In this article, we first consider the case when $\Delta$ has only one totally interior edge, because it is the simplest non-trivial case. Then we may apply the formula we find here to get an upper bound on some more general cases.
\end{abstract}

\section{Introduction}
\subsection{Spline modules}
Let $\Delta$ be a connected, pure $2$-dimensional simplicial complex embedded in $\RR^2$. Throughout this article, we assume that the genus of $\Delta$ is $0$. Let $r\geq0$ be an integer. We define $C^{r}(\Delta)$ to be the set of $C^{r}$-differentiable piecewise polynomial functions on $\Delta$. These functions are called splines. More explicitly:
\begin{Definition}
$\SMr$ is the set of functions $f:\Delta\rightarrow \RR$ such that:
\begin{enumerate}
\item For all facets $\sigma\in\Delta$, $f|_{\sigma}$ is a polynomial in $\RR[x,y]$.
\item $f$ is differentiable of order $r$.
\end{enumerate}
\end{Definition}
Therefore, $\SMr$ is a module over $\RR[x,y]$. For each integer $d\geq 0$, we define
\begin{align*}
C^{r}_{d}(\Delta):=\{f\in C^{r}(\Delta):\deg(f|_\sigma)\leq d, \text{ for all facets }\sigma\in\Delta\}.
\end{align*}
It is a finite dimensional $\RR$-vector space. One of the key problems in spline theory is the determination of the dimensions of $C^{r}_{d}(\Delta)$ for all $d$.\\
  Assume that $S=\RR[x,y,z]$. Let $\hat{\Delta}$ be the cone of $\Delta$, i.e., embed $\Delta$ in $\mathbb{R}^3$ by sending $(x,y)$ to $(x,y,1)$ and define $\hat{\Delta}$ to be the cone of $\Delta$ over the origin of $\RR^3$. Then $C^{r}(\hat{\Delta})$, the so-called homogenized spline module of $\Delta$, is a graded $S$-module and  $\dim_\RR C^{r}_{d}(\Delta)=\dim_\RR C^{r}(\hat{\Delta})_d=\HF(C^{r}(\hat{\Delta}),d)$, where $\HF$ stands for the Hilbert function of a graded $S$-module.
\subsection{The homological algebra tools}
Billera introduced the use of homological algebra in spline theory in \cite{billera1988homology}. Following this path, Schenck and Stillman\cite{schenck1997local} defined a chain complex $J_{\bullet}$ to deal with the problem of freeness of $\SMr$. Consider $S_{\bullet}$ as the relative simplicial complex $C_\bullet(\Delta,\partial\Delta)$ with coefficients in $S$. 
Denote by $l_{\hat{\tau}}$ a linear form vanishing on $\hat{\tau}\in\Delta_1$. The authors of \cite{schenck1997local} define the ideal complex $J_\bullet$ to be
\begin{equation}
\label{idealComplex}
J_\bullet: 0\rightarrow \bigoplus_{\tau\in\Delta^{\circ}_1} J(\tau)\xrightarrow{\partial} \bigoplus_{v\in\Delta^{\circ}_0}J(v)\rightarrow 0,
\end{equation}
where
\begin{align*}
J(\tau)=\langle l_{\hat{\tau}}\rangle^{r+1}\\
J(v)=\sum_{v\in\tau} J(\tau)
\end{align*}
are ideals of $S$, 
and $\partial$ is induced by the differential in $S_{\bullet}$. So $J_{\bullet}$ is a subcomplex of $S_{\bullet}$, and we may consider the quotient complex $S_\bullet/J_\bullet$. It has a lot of good features. First of all, as an $S$-module, $\hSMr$ is isomorphic to $H_2(S_\bullet/J_\bullet)$. Second the Hilbert function of $\SMr$ may be obtained from local data of $\Delta$ and the Hilbert function of $H_1(S_\bullet/J_\bullet)$.
So the problem is to calculate $\HF(H_1(S_\bullet/J_\bullet),d)$. In particular, since we know that the length of $H_1(S_\bullet/J_\bullet)$ is finite, we'd like to know the least $d$ such that $(H_1(S_\bullet/J_\bullet))_d=0$. More precisely, we ask for an estimate for the Castelnuovo-Mumford regularity of $H_1(S_\bullet/J_\bullet)$.
\begin{Definition}[Castelnuovo-Mumford regularity]
	Assume $M$ is an $S$-module of finite length. Then the Castelnuovo-Mumford regularity may be defined as
	\begin{equation*}
	\reg M:=\max\{M_d\neq 0, d\geq 0\}.
	\end{equation*}
\end{Definition}

\subsection{Known upper bounds and the ``$2r+1$" conjecture on regularity of $H_0(J_\bullet)$}
In \cite{alfeld1987dimension}, Alfeld and Schumaker show that the regularity of $H_1(S_\bullet/J_\bullet)$ is less than $4r+1$. They improve the result to $3r+1$ in a later paper\cite{alfeld1990dimension}. Schenck and Stiller conjectured in \cite{schenck2002cohomology}:
\begin{Conjecture}[Schenck-Stiller]
	\begin{equation*}
		\reg H_1(S_{\bullet}/J_{\bullet})\leq 2r.
	\end{equation*}
\end{Conjecture}
They call it the ``$2r+1$" conjecture, because there is an equivalent statement saying that $H_1(S_{\bullet}/J_{\bullet})$ vanishes at degrees greater than or equal to $2r+1$. 
Work of Toh\v{a}neanu\cite{tohaneanu2005smooth}  shows that this guess is optimal by finding a $\Delta$ such that $\reg  H_1(S_{\bullet}/J_{\bullet})=2r$.
\begin{figure}[h]
		\centering
		\includegraphics[width=0.5\textwidth]{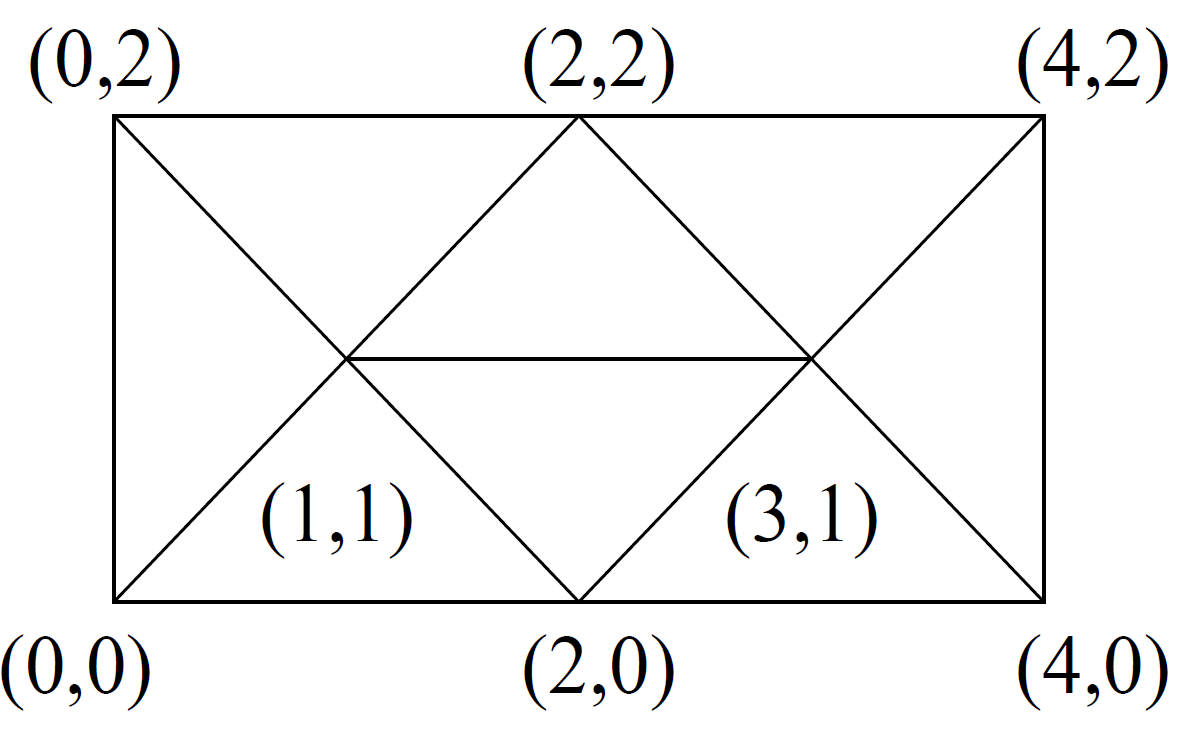}
		\caption{Toh\v{a}neanu's $\Delta$ such that $\reg  H_1(S_{\bullet}/J_{\bullet})=2r$.\cite{tohaneanu2005smooth}}
		\label{TohConfig}
\end{figure}

Note that in our case, the genus of $\Delta$ is $0$, it is always true that $H_0(J_\bullet)=H_1(S_\bullet/J_\bullet)$. So we will alway compute $H_0(J_\bullet)$, instead of $H_1(S_\bullet/J_\bullet)$.\\

In order to state our results, first we need to introduce some notions. We call an edge $\varepsilon=[v_1v_2]$ \textbf{totally interior edge} if both $v_1$ and $v_2$ are interior vertices of $\Delta$. If either $v_1$ or $v_2$ is interior but not both, then we call $\varepsilon$ a \textbf{partially interior edge}. We also adopt the notations from \cite{schenck2002cohomology}: For an interior vertex $v\in\Delta^{\circ}_0$, let $f_1(v)$ denote the number of edges incident to $v$, $k(v)$ the number of those edges of distinct slope. Let $\numtotalintegde(v)$ denote the number of totally interior edges incident to $v$, and $k^{00}(v)$ be the number of those edges of distinct slope. Let $\numpartialedge(v)$ be the number of edges incident to $v$ with one boundary vertex and one interior vertex, and $k^{0\partial}(v)$ be the number of those edges of distinct slope. And we assume that if $\varepsilon_1\in\Delta_1^{0\partial}(v)$ and $\varepsilon_0\in\Delta_1^{00}(v)$, then they have different slopes. Let $\alpha(v)=\lfloor\frac{r+1}{k^{0\partial}(v)}\rfloor$.\\
\begin{Example}
\begin{figure}[h]
		\centering
		\includegraphics[width=0.5\textwidth]{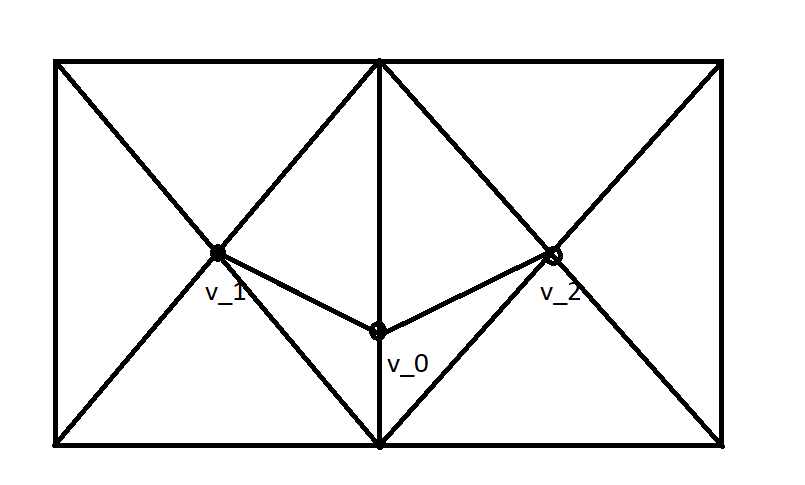}
		\caption{A configuration with $2$ totally interior edges}
		\label{ce1}
\end{figure}
For the configuration in figure \ref{ce1}, the numbers of totally interior edges incident to a vertex are $f_1^{00}(v_1)=f_1^{00}(v_2)=1$ and $f_1^{00}(v_0)=2$. The number of partial interior edges $\numpartialedge(v_1)=4$, but the number of slopes $k^{0\partial}(v)=2$ is different from $\numpartialedge(v_1)$. The number of slopes $k(v_1)=3$.\\
In this example, $\alpha(v_1)=\alpha(v_2)=\lfloor\frac{r+1}{2}\rfloor$.
\end{Example}
With these notations, the main theorem of this article may be state as following:
\begin{Theorem}[Main theorem]
	If $\Delta$ has only one total interior edge $\varepsilon=[v_1v_2]$, then
	\begin{equation}
	\alpha(v_1)+\alpha(v_2)+r-1 \leq \reg H_0(J_\bullet)\leq \alpha(v_1)+\alpha(v_2)+r.
	\end{equation}
\end{Theorem}
As applications, we also show that
\begin{Theorem}
	\label{Thm1tot}
	Let $\Delta$ be a simplicial complex with genus $0$. If $\Delta$ has only one totally interior edge, then $\reg H_0(J_\bullet)\leq 2r$.
\end{Theorem}
and 
\begin{Theorem}\label{ThmPath}
	If $\numpartialedge(v)\neq 0$ for each $v\in\intv$, then
	\begin{equation*}
		\max_{\varepsilon=[v_iv_j]\in\totintegde}\{\lfloor\frac{r+1}{k(v_i)-1}\rfloor+\lfloor\frac{r+1}{k(v_j)-1}\rfloor+r-1\}\leq \reg H_0(J_\bullet)\leq  \max_{\varepsilon=[v_iv_j]\in\totintegde}\{\alpha(v_i)+\alpha(v_j)+r\}.
	\end{equation*}
\end{Theorem}

\section{Proof of the main theorem}
Without loss of generality, we only need to consider the case when $\Delta$ is \textbf{indecomposable}, that is, when the first differential map $\partial_1$ of $S_\bullet$ 
\begin{equation*}
	\partial_1:\bigoplus_{\tau\in\Delta^{\circ}_1} S\rightarrow\bigoplus_{v\in\Delta^{\circ}_0} S
\end{equation*}
is indecomposable. This is because if $\Delta$ is decomposable, we may just consider all its indecomposable components. From now on, if not otherwise claim, we will assume all $\Delta$ we consider is indecomposable.

\subsection{Results on $\Delta=\text{star}(v)$}
If there is only one interior vertex $v$ in $\Delta$, we will say $\Delta$ is a star of $v$ and denote $\Delta=\text{star}(v)$. After a change of coordinates, we may move $v$ to the origin. Then $J(v)$ is an ideal in $2$ variables, and since each vertex we have at least two edges with different slopes, so $S/J(v)$ has projective dimension $2$. In \cite{schumaker1979dimension}, Schmaker gave a dimension formula for the star, from which it follows that
\begin{Theorem}[Schumaker]
	\label{Thm_LocalData}
	A free resolution of $S/J(v)$ is given by
	\begin{equation}
		\label{SESofLocalData}
		S(-r-1-\alpha(v))^{a_1}\oplus S(-r-2-\alpha(v))^{a_2}\rightarrow S(-r-1)^{k(v)}\rightarrow S\rightarrow S/J(v)\rightarrow 0,
	\end{equation}
	where $\alpha(v)=\lfloor (r+1)/(k(v)-1)\rfloor$, $a_1=(k(v)-1)\alpha(v)+k(v)-r-2$ and $a_2=k(v)-1-a_1=r+1-(k(v)-1)\alpha(v)$.
\end{Theorem}

\subsection{A presentation of $H_{0}(J_{\bullet})$}
Let
\begin{align*}
	&F_{\bullet}: \oplus S(-r-1),\\
	&G_{\bullet}: \oplus S(-r-1-m)^{a_1}\oplus S(-r-2-m)^{a_2}\rightarrow \oplus S(-r-1)^{a}.
\end{align*}
be free resolutions of $\oplus J(\tau)$ and $\oplus J(v)$, respectively. Then the map $\partial :\oplus J(\tau)\rightarrow \oplus J(v)$ induces a map between their free resolutions:
\begin{equation}
	\mcdiff : F_{\bullet}\rightarrow G_{\bullet}.
\end{equation}
The mapping cone $P_{\bullet}$ of $\mcdiff$ has
\begin{equation}
	H_0(P_{\bullet})=H_{0}(J_{\bullet}).
\end{equation}

\subsection{The $1$-total-interior-edge case}
If there is only one total interior edge $\varepsilon=[v_1v_2]$ in $\Delta$, then $v_1$ and $v_2$ are the only interior vertices in $\Delta$. To simplify the notation, we assume the number of distinct slopes $k(v_1)=a$, $k(v_2)=b$ and $a\leq b$ in this case.
\begin{Lemma}
	If $\Delta$ has only one total interior edge $\varepsilon=[v_1v_2]$, let
	\begin{align*}
		J'(v_1)=\sum_{v_1\in\tau,\tau\neq \varepsilon}J(\tau),\\
		J'(v_2)=\sum_{v_2\in\tau,\tau\neq \varepsilon}J(\tau).
	\end{align*} 
	Then
	\begin{equation}
		H_0(J_\bullet)\simeq S/(J'(v_1):J(\varepsilon)+J'(v_2):J(\varepsilon)).
	\end{equation}
\end{Lemma}
\begin{proof}
	By choosing coordinates, we may assume that 
	\begin{equation}
		(G_2\to G_1)=\begin{bmatrix}
		(g_{ij})_{a\times(a-1)}& 0\\
		0& (h_{ij})_{b\times(b-1)}
		\end{bmatrix},
	\end{equation}
	hence
	\begin{equation}
		(P_1\to P_0)=\begin{bmatrix}
		I_{(a-1)\times(a-1)}& &0&(g_{ij})_{1\leq i,j\leq a-1}&0\\
		0&1&0&(g_{aj})_{1\leq j\leq a-1} &0\\
		0&1&0&0 &(h_{1\leq j\leq b-1})\\
		0& &I_{(b-1)\times(b-1)}&0 &(h_{ij})_{1\leq j\leq b-1,2\leq i\leq b}
		\end{bmatrix}.
	\end{equation}
	Therefore, $\text{coker}(P_1\rightarrow P_0)$ is isomorphic of the cokernel of
	\begin{equation}
		\begin{bmatrix}
		(g_{aj})_{1\leq j\leq a-1} &(h_{1\leq j\leq b-1})
		\end{bmatrix}.
	\end{equation}
	The image of $\begin{bmatrix}
	(g_{aj})_{1\leq j\leq a-1}
	\end{bmatrix}$ is $J'(v_1):J(\varepsilon)$ and that of $\begin{bmatrix}
	(h_{1\leq j\leq b-1})
	\end{bmatrix}$ is $J'(v_2):J(\varepsilon)$. Therefore,
	\begin{equation}
		\text{coker}(P_1\rightarrow P_0)\simeq S/(J'(v_1):J(\varepsilon)+J'(v_2):J(\varepsilon)).
	\end{equation}
\end{proof}

\begin{Definition}[Initial ideal]
After choosing a monomial order of the ring $\RR[x,y,z]$, we denote by $\text{in}(g)$ the leading term of a polynomial $g$. The initial ideal of $I$ is defined as
\begin{equation*}
\text{In}(I)=\{\text{in}(g)|g\in I\}.
\end{equation*}
\end{Definition}
The initial ideal has many good properties. It is a monomial ideal with the same Hilbert function as that of $I$.
\begin{Lemma}\label{lemma_fullrank}
Let 
\begin{equation*}
J=\langle (x+{c_1}y)^{r+1},(x+{c_2}y)^{r+1},\dots,(x+{c_s}y)^{r+1}\rangle
\end{equation*}
be an ideal in $R=\RR[x,y]$, where $c_1,\dots,c_s$ are distinct constants in $\RR$ and $s\geq 2$. Then with respect to the lexicographic order,
\begin{align}\label{InJ}
	\text{In}(J)=\langle x^{r+1},&x^{r}y,\dots, x^{r-(s-1)+1}y^{s-1},\nonumber\\ &x^{r-(s-1)}y^{s+1},\dots,x^{r-2(s-1)+1}y^{2s-1},\nonumber\\
	&\dots\\
	&x^{r-j(s-1)}y^{js+1},\dots,x^{r-(j+1)(s-1)+1}y^{(j+1)s-1},\nonumber\\
	&\dots\nonumber\\
	&x^{r-\lfloor\frac{r}{s-1}\rfloor(s-1)}y^{\lfloor\frac{r}{s-1}\rfloor s+1},\dots,y^{r+1+\lfloor\frac{r}{s-1}\rfloor}\rangle.\nonumber
	\end{align}
\end{Lemma} 
\begin{proof}
We only consider $d\geq r+1$. Take the basis
\begin{equation*}
(x^{d},x^{d-1}y,\dots,xy^{d-1},y^{d})
\end{equation*}
for $R_{d}$. Then $J_d$ is the image of the matrix
\begin{equation*}
A=\begin{bmatrix}
A_1&A_2&\dots&A_s
\end{bmatrix},
\end{equation*}
where $A_i$ is a $(d+1)\times (d-r)$ matrix
\begin{equation*}
A_i=\begin{bmatrix}
1&&&\\
\binom{d}{1}c_i&1&&\\
\binom{d}{2}c_i^{2}&\binom{d-1}{1}c_i&&1\\
\dots&&&\\
\binom{d}{d}c_i^{d}&\binom{d-1}{d-1}c_i^{d-1}&\dots&\binom{r+1}{r+1}c_i^{r+1}
\end{bmatrix}~\text{for}~i=1,2,\dots,s.
\end{equation*}
By the proof of Theorem 2.1 in \cite{schumaker1979dimension}, 
\begin{equation*}
\text{rank}A=\min\{d+1,s\cdot (d-r)\}.
\end{equation*}
In other words, the matrix $A$ is always of full rank.\\
Therefore, for $d\geq \frac{sr+1}{s-1}$, i.e. $d+1\geq s\cdot(d-r)$, $J_d=S_d$.\\
Now we consider $d<\frac{sr+1}{s-1}$. In this case, we want to show that the top-most maximal minor of $A$ is non-zero. If $A'$ is the submatrix of $A$ formed by the first $s\cdot(d-r)$ rows of $A$. Explicitly,
\begin{equation*}
A'=\begin{bmatrix}
A_1'&A_2'&\dots&A_s'
\end{bmatrix},
\end{equation*}
where $A_i'$ is a $s(d-r)\times (d-r)$ matrix
\begin{equation*}
A_i=\begin{bmatrix}
1&&&\\
\binom{d}{1}c_i&1&&\\
\binom{d}{2}c_i^{2}&\binom{d-1}{1}c_i&&1\\
\dots&&&\\
\binom{d}{s(d-r)-1}c_i^{s(d-r)-1}&\binom{d-1}{s(d-r)-2}c_i^{s(d-r)-2}&\dots&\binom{r+1}{s(d-r)-d+r}c_i^{s(d-r)-d+r}
\end{bmatrix}
\end{equation*}
for $i=1,2,\dots,s$.\\
By the same reasoning as in \cite{schumaker1979dimension}, except for the first column, each column is obtained by differentiating its predecessor. So $A'$ is in fact the matrix corresponding to Hermite interpolation at the points $c_1,\dots,c_s$ with respect to
\begin{equation*}
1,\binom{d}{1}t,\binom{d}{2}t^2,\dots,\binom{d}{s(d-r)}t^{s(d-r)-1}.
\end{equation*}
These are just scalar multiples of power functions. So 
\begin{equation*}
\text{rank}A'= s(d-r).
\end{equation*}
This means the leading terms of elements in $J_{d}$ are
\begin{equation*}
	\{x^{d-j}y^{j}:0\leq j\leq s(d-r)-1\}.
\end{equation*}
So $\text{In}(J)$ is of the form ($\ref{InJ}$). 
In particular, $\text{In}(J)$ is a lex-segment ideal. 
\end{proof}

\begin{Lemma}
	By choosing coordinates, we may assume that $v_1=[0,1,0]$ and $v_2=[1,0,0]$, and therefore
	\begin{align*}
		J(\varepsilon)&=\langle z^{r+1}\rangle,\\
		J'(v_1)&=\langle (x+{\alpha_0}z)^{r+1},(x+{\alpha_1}z)^{r+1},\dots,(x+{\alpha_{(a-2)}}z)^{r+1}\rangle,\\
		J'(v_2)&=\langle (y+{\beta_0}z)^{r+1},(y+{\beta_1}z)^{r+1},\dots,(y+\beta_{(b-2)}z)^{r+1}\rangle,
	\end{align*}
	where we can always choose $\alpha_0=\beta_0=0$. 
	\begin{equation}
	\label{eq_initial}
		\text{In}(J'(v_1):J(\varepsilon)+J'(v_2):J(\varepsilon))=\text{In}(J'(v_1):J(\varepsilon))+\text{In}(J'(v_2):J(\varepsilon)),
	\end{equation}
	and
	\begin{equation}
	\label{eq_initial_quotient}
		\text{In}(J'(v_1):J(\varepsilon))=\text{In}(J'(v_1)):J(\varepsilon).
	\end{equation}
\end{Lemma}
\begin{proof}
	By Lemma \ref{lemma_fullrank}, with respect to the lexicographic order,
	\begin{align*}
	\text{In}(J'(v_1))=\langle x^{r+1},&x^{r}z,\dots, x^{r-(a-2)+1}z^{a-2},\\ &x^{r-(a-2)}z^{a},\dots,x^{r-2(a-2)+1}z^{2a-3},\\
	&\dots\\
	&x^{r-j(a-2)}z^{j(a-1)+1},\dots,x^{r-(j+1)(a-2)+1}z^{(j+1)(a-1)-1},\\
	&\dots\\
	&x^{r-\lfloor\frac{r}{a-2}\rfloor(a-2)}z^{\lfloor\frac{r}{a-2}\rfloor (a-1)+1},\dots,z^{r+1+\lfloor\frac{r}{a-2}\rfloor}\rangle.
	\end{align*}
	Let $n=\lfloor\frac{r+1}{a-1}\rfloor$, then
	\begin{align*}
	\text{In}(J'(v_1):(z^{r+1}))
	=\langle &x^{n},\dots,x^{r-(n+1)(a-2)+1}z^{(n+1)(a-1)-r-2},\\
	&x^{r-(n+1)(a-2)}z^{(n+1)(a-1)-r},\dots,x^{r-(n+2)(a-2)+1}z^{(n+2)(a-1)-r-2},\\
	&\dots\\
	&x^{r-\lfloor\frac{r}{a-2}\rfloor(a-2)}z^{\lfloor\frac{r}{a-2}\rfloor (a-1)-r},\dots,z^{\lfloor\frac{r}{a-2}\rfloor}\rangle\\
	&=\text{In}(J'(v_1)):(z^{r+1}).
	\end{align*}
	This proves (\ref{eq_initial_quotient}).\\
	Following Gruson and Peskine's convention\cite{gruson1978genre}, we denote
	\begin{equation*}
	\text{In}(J'(v_1))=\langle x^{r+1},x^{r}z^{\lambda_r},x^{r-1}z^{\lambda_{r-1}},\dots,xz^{\lambda_1},z^{\lambda_0}\rangle,
	\end{equation*}
	and
	\begin{equation*}
	\text{In}(J'(v_2))=\langle y^{r+1},y^{r}z^{\eta_r},y^{r-1}z^{\eta_{r-1}},\dots,yz^{\eta_1},z^{\eta_0}\rangle,
	\end{equation*}
	and let
	\begin{align*}
	\upsilon_i:=\lambda_{r-i}-\lambda_{r-i+1},\\
	\epsilon_i:=\eta_{r-i}-\eta_{r-i+1}.
	\end{align*}
	Then
	\begin{equation*}
	\lambda_{r}=\eta_{r}=1,
	\end{equation*}
	and
	\begin{align*}
	\upsilon_i=\begin{cases}
	2,\text{ if }i\equiv 0\mod a-2,\\
	1,\text{ otherwise.}
	\end{cases}\\
	\epsilon_i=\begin{cases}
	2,\text{ if }i\equiv 0\mod b-2,\\
	1,\text{ otherwise.}
	\end{cases}
	\end{align*}
	Let
	\begin{align*}
	\lambda'_{i}:=\max\{0,\lambda_i-r-1\},\\
	\eta'_{i}:=\max\{0,\eta_i-r-1\},
	\end{align*}
	and $i_0$( $j_0$ resp.) be the least $i$ such that $\lambda_i\leq r+1$( $\eta_i>r+1$ resp.), so
	\begin{align*}
	i_0=\lfloor\frac{r+1}{a-1}\rfloor,\\
	j_0=\lfloor\frac{r+1}{b-1}\rfloor.
	\end{align*} 
	Then
	\begin{align*}
	\text{In}(J'(v_1):(z^{r+1}))=\langle x^{i_0}z^{\lambda'_{i_{0}}},x^{i_0-1}z^{\lambda'_{i_{0}-1}},\dots,xz^{\lambda'_{1}},z^{\lambda'_{0}}\rangle,\\
	\text{In}(J'(v_2):(z^{r+1}))=\langle y^{j_0}z^{\eta'_{j_{0}}},y^{j_0-1}z^{\eta'_{j_{0}-1}},\dots,yz^{\eta'_{1}},z^{\eta'_{0}}\rangle,
	\end{align*}
	where $i_0=j_0=0$. 
	Notice that for $0\leq i<r-i_0$,
	\begin{equation*}
	\lambda'_{i}-\lambda'_{i+1}=\lambda_{i}-\lambda_{i+1}=\upsilon_{r-i}.
	\end{equation*}
	Denote $J'(v_1):J(\varepsilon)$ and $J'(v_2):J(\varepsilon)$ by $Q(v_1)$ and $Q(v_2)$, respectively.\\
	Now let $g\in Q(v_1)$ and $h\in Q(v_2)$ be homogeneous elements, then $g$ is in $x$ and $z$ and $h$ is in $y$ and $z$. So $\text{gcd}(f,g)=\text{gcd}(\text{In}(f),\text{In}g)$, and hence
	\begin{equation*}
		\text{gcd}(\text{In}\frac{g}{\text{gcd}(g,h)},\text{In}\frac{h}{\text{gcd}(g,h)})=1.
	\end{equation*}
	This means that the $S$-pair of $g$ and $h$ reduces to $0$. This proves (\ref{eq_initial}).\\
	Let $Q=Q(v_1)+Q(v_2)$ and $l_0$ be the smallest $i$ such that $\lambda'_{i}<\eta'_0$. Then
	\begin{align*}
	\text{In}Q=\langle &x^{i_0},x^{i_0-1}z^{\lambda'_{i_0-1}},\dots,x^{l_0}z^{\lambda'_{l_0}},\\
	&y^{j_0},y^{j_0-1}z^{\eta'_{j_0-1}},\dots,yz^{\eta'_1},z^{\eta'_0}\rangle.
	\end{align*}
\end{proof}
\begin{Remark}
	Here is an estimation of $l_0$: The inequality
	\begin{equation*}
		\lfloor\frac{r-l_0}{a-2}\rfloor-l_0-r-1=\lambda'_{l_0}<\eta'_0=\lfloor\frac{r}{b-2}\rfloor-r-1
	\end{equation*}
	implies
	\begin{equation*}
	\frac{r-l_0}{a-2}-\frac{a-3}{a-2}-l_0\leq\frac{r}{b-2}.
	\end{equation*}
	Therefore we get a lower bound,
	\begin{equation*}
		l_0>\frac{b-a}{(a-1)(b-2)}r-\frac{(a-3)}{(a-1)}.
	\end{equation*}
\end{Remark}
Next, we want to describe the higher syzygies of $\text{In}Q$. First, let's illustrate how they look like with an example. Recall from \cite{GTM227} that a \textbf{Buchberger graph} $\text{Buch}(I)$ of a monomial ideal $I=\langle m_1,\dots,m_n\rangle$ has vertices $i=1,\dots, n$ assigned with $m_i$ and an edge $(i,j)$ whenever there is no monomial $m_k$ divides $\text{lcm}(m_i,m_j)$. The $k$-face $(i_1,\dots,i_k)$ of the graph is assigned with the monomial $\text{lcm}(m_{i_1},\dots,m_{i_k})$.
\begin{Example}
	\begin{figure}
		\centering
		\includegraphics[width=0.5\textwidth]{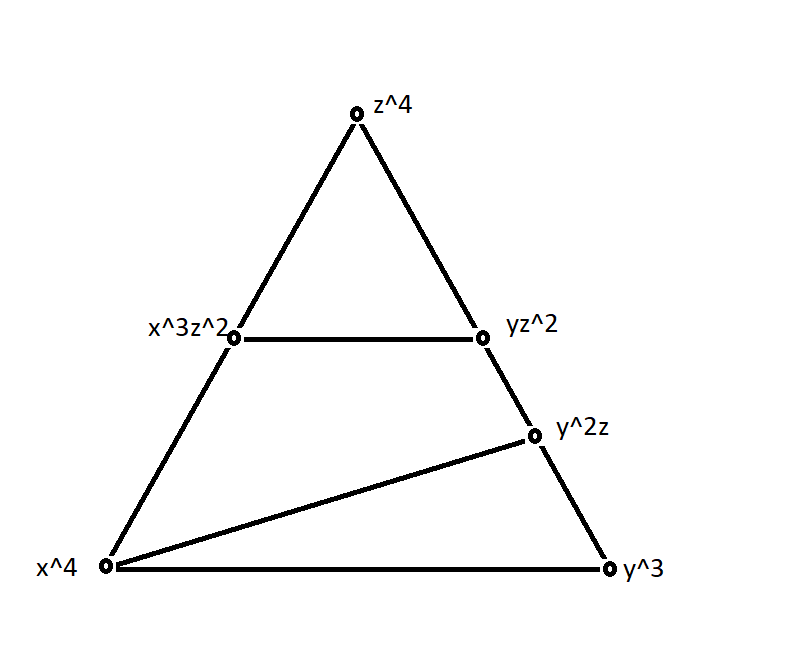}
		\caption{Buchberger graph of $\text{In}Q$ when $(a,b)=(3,4)$ and $r=8$}
		\label{a3b4r08}
	\end{figure}
Assume that $(a,b)=(3,4)$ and $r=8$. Then $i_0=4$ and $j_0=3$, 
\begin{align*}
\text{In}(J'(v_1))=\langle x^9, x^8z, x^7z^3,x^6z^5,x^5z^7,x^4z^9,x^3z^{11},x^2z^{15},z^{17}\rangle
\end{align*}
and
\begin{equation*}
\text{In}(J'(v_2))=\langle y^9, y^8z, y^7z^2,y^6z^4,y^5z^5,y^4z^7,y^3z^8,y^2z^{10},yz^{11},z^{13}\rangle.
\end{equation*}
Therefore,
\begin{equation*}
\text{In}Q=\langle x^4, x^3z^2,y^3,y^2z,yz^2,z^4\rangle.
\end{equation*}
Its Buchberger graph is in Figure \ref{a3b4r08}. We can read all its higher syzygies from the graph. The second syzygy is generated by the edges, 
\begin{equation*}
	x^4z^2,x^3z^4,y^3z,y^2z^2,yz^4,x^4y^3,x^4y^2z,x^3yz^2.
\end{equation*}
The third syzygy is generated by the faces,
\begin{equation*}
	x^4y^3z,x^4y^2z^2,x^3yz^4.
\end{equation*}
The regularity of $\text{In}Q$ is measured at the third syzygy.
\begin{equation*}
	\reg \text{In}Q=\deg(x^4y^3z)-3=5.
\end{equation*}
\end{Example}
What we observed from the example can be generalized as following:
\begin{Lemma}
	The second syzygy of $\text{In}Q$ is generated by the union of the following set:
	\begin{enumerate}
		\item $\{x^iz^{\lambda'_{i-1}}: l_0< i\leq i_0\}\cup\{x^{l_0}z^{\eta'_0}\}$;
		\item $\{y^jz^{\eta'_{j-1}}: 1\leq j\leq j_0\}$;
		\item $\{x^iy^jz^{\eta'_j}: \text{ for all } (i,j) \text{ such that }\lambda'_i\leq\eta'_j<\lambda'_{i-1}\}$;
		\item $\{x^iy^jz^{\lambda'_i}: \text{ for all } (i,j) \text{ such that }\eta'_j\leq\lambda'_i<\eta'_{j-1}\}$.
	\end{enumerate}
\end{Lemma}
\begin{proof}
	Each of the elements listed above is an l.c.m. of two elements of the first syzygy. And it cannot be divided by other elements of the first syzygy. We only need to show that they generate the second syzygy.\\
	If $x^lz^t\in\text{Syz}_2$, then it must be l.c.m. of two adjacent elements, i.e. $x^{i}z^{\lambda'_{i}}$ and $x^{i-1}\lambda'_{i-1}$. The only exception is the l.c.m. of $x^{l_0}z^{\lambda_{l_0}}$ and $z^{\eta'_0}$. Similar reasoning applies to $y^lz^t\in\text{Syz}_2$. They form the list 1. and 2.\\
	Now we consider elements $x^iy^jz^t\in\text{Syz}_2$ with $i,j\neq 0$. If $\eta'_{j}\geq \lambda'_{i-1}$, then the l.c.m. of $x^iz^{\lambda'_i}$ and $y^jz^{\eta'_j}$ can be divided by $x^{i-1}z^{\lambda'_{i-1}}$. Similar reasoning applies to $\lambda'_i\leq\eta'_{j-1}$. So elements of the form $x^iy^jz^t$ with nonzero $i$ and $j$ must fall in either list 3. or 4.\\
	So the above list is a full list of elements in $\text{Syz}_2$.
\end{proof}
\begin{Corollary}
	The Buchberger graph of $\text{In}Q$ is planar.
\end{Corollary}
\begin{Lemma}\label{lemma7}
	Any two monomial 
	generators of the third syzygy of $\text{In}Q$ are in different degree of $z$. If we order them by the degree of $z$:
	\begin{equation*}
		h_0(x,y)z^{\zeta_0}>h_1(x,y)z^{\zeta_1}>\dots>h_{t}(x,y)z^{\zeta_t},
	\end{equation*}
	where $\zeta_0<\zeta_1<\dots<\zeta_t$, then the total degree has the following property:
	\begin{equation*}
		\deg h_0(x,y)z^{\zeta_0}\geq \deg h_1(x,y)z^{\zeta_1}\geq\dots\geq \deg h(x,y)_{t}z^{\zeta_t}.
	\end{equation*}
	So
	\begin{equation*}
		\reg  H_0(J_\bullet)=\deg h_0(x,y)z^{\zeta_0}-3+(r+1).
	\end{equation*}
	In other words, the regularity of $H_0(J_\bullet)$ is always measured at the "bottom" face of its Buchberger graph.
\end{Lemma}
\begin{proof}
We only need to prove the first statement. In fact, the degree of $z$ in the monomial assigned to a face is determined by the highest vertex of that face. So we may order the faces by the height of thier highest vertex. The order is strict, because if we have $\zeta=\zeta'=t$, the situation can only be the one in figure \ref{bgnonexist}.\\
\begin{figure}[h]
		\centering
		\includegraphics[width=0.5\textwidth]{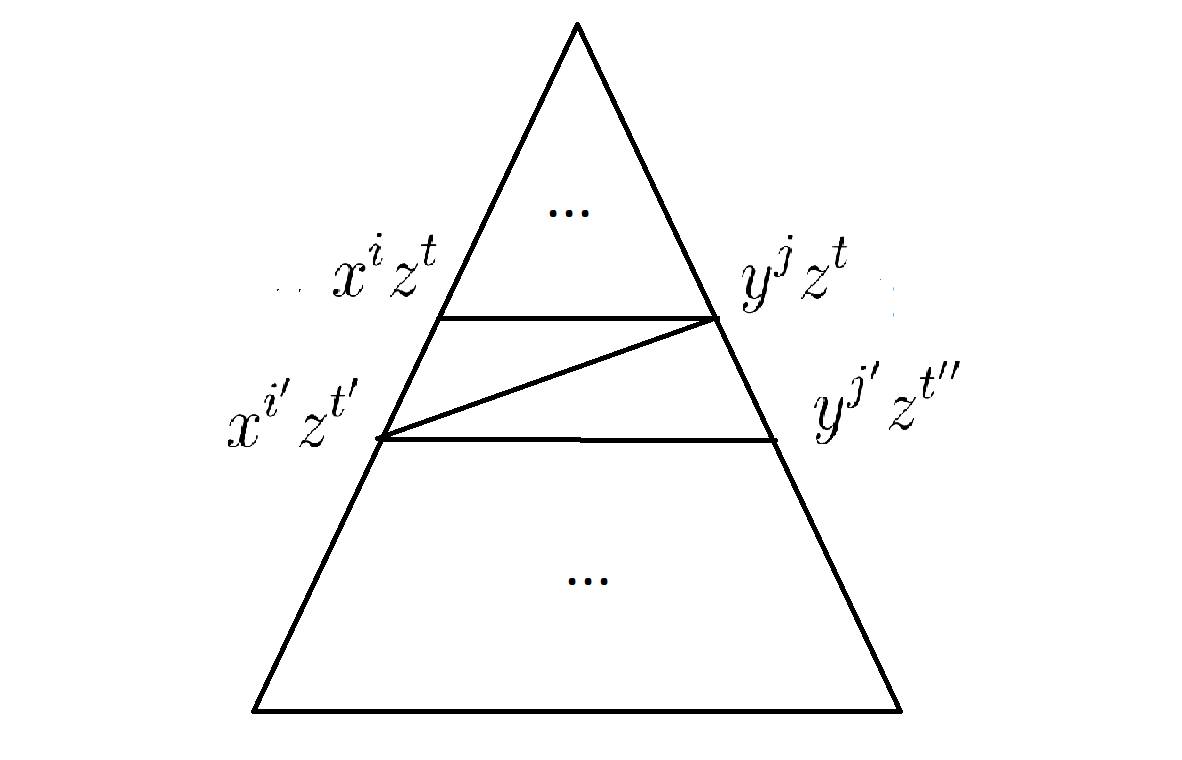}
		\caption{This type of Buchberger graph $\text{Buch}(Q)$ does not exist}
		\label{bgnonexist}
\end{figure}
Notice that $i'>i$. This type of $Buch(Q)$ cannot exist, because $x^iz^t$ divides $\text{lcm}(x^{i'}z^{t'},y^jz^{t})$, which is a contradiction to the definition of the Buchberger graph.
\end{proof}

\begin{Proposition}
Let $(a,b)$ be the number of slopes at two interior vertices (including the totally interior one). The regularity of $H_0(J_{\bullet})$ is bounded by
\begin{equation}
 \lfloor\frac{r+1}{a-1}\rfloor+\lfloor\frac{r+1}{b-1}\rfloor+r-1\leq \reg H_0(J_{\bullet})\leq\lfloor\frac{r+1}{a-1}\rfloor+\lfloor\frac{r+1}{b-1}\rfloor+r.
\end{equation}
\end{Proposition}
\begin{proof}
	Because $H_0(J_\bullet)$ is Artinian, so the regularity is measured at its last syzygy (for a more general statement, c.f. \cite{GTM229}, 4E.5). By Lemma \ref{lemma7}, 
	\begin{equation*}
	\reg  H_0(J_\bullet)=\deg h_0(x,y)z^{\zeta_0}-3+(r+1),
	\end{equation*}
	where $h_0(x,y)=x^{i_0}y^{j_0}$ and $\zeta_0=\min\{\lambda'_{i_0-1},\eta'_{j_0-1}\}\in\{1,2\}$.
\end{proof}

\section{Proof of Theorem \ref{Thm1tot}}
We may apply the main theorem to prove that, for $a\geq 3$ and $b\geq 4$, Theorem \ref{Thm1tot} holds, because we can see that 
\begin{equation*}
\lfloor\frac{r+1}{2}\rfloor+\lfloor\frac{r+1}{3}\rfloor+r\leq 2r,
\end{equation*}
for all positive integer $r$. We know that $b\geq a\geq 3$, because the cases we consider are triangularization. So we only need to check the case when $(a,b)=(3,3)$.
\begin{figure}
\centering
\includegraphics[width=0.5\textwidth]{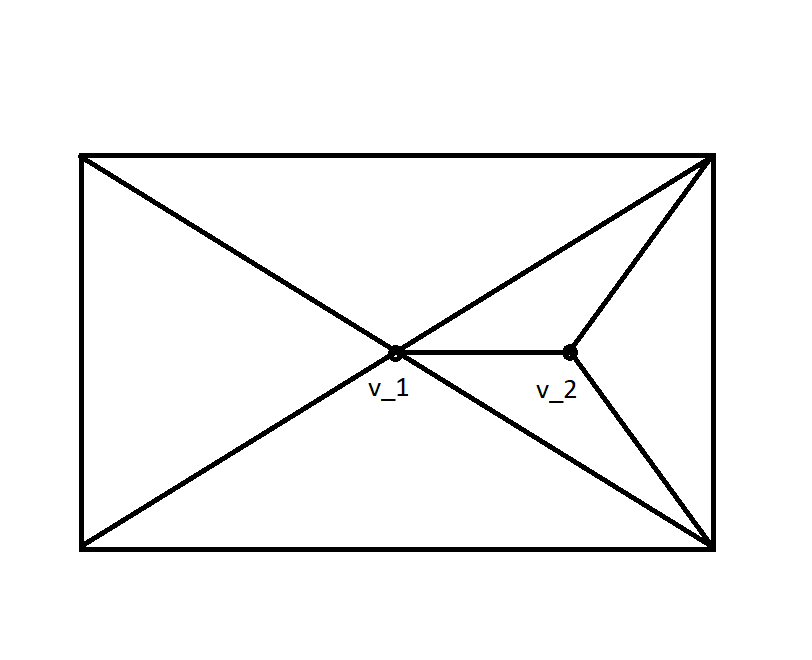}
\caption{Case $(a,b)=(3,3)$}
\label{fig_3_3}
\end{figure}
\begin{Proposition}
Assume that $(a,b)=(3,3)$. Then
\begin{equation*}
\reg H_0(J_{\bullet})=2r.
\end{equation*}
\end{Proposition}
\begin{proof}
\begin{itemize}
\item If $r$ is even, then let $r=2p$.
\begin{align*}
\text{In}Q(v_1)=\langle x^p,x^{p-1}z^2,\dots,xz^{2p-2},z^{2p}\rangle\\
\text{In}Q(v_2)=\langle y^p,y^{p-1}z^2,\dots,yz^{2p-2},z^{2p}\rangle
\end{align*}
By lemma \ref{lemma7}, 
\begin{equation*}
\reg  H_0(J_\bullet)=\deg h_0(x,y)z^{\zeta_0}-3+(r+1)=2p+2-3+2p+1=4p=2r.
\end{equation*}
\item If $r$ is odd, then let $r=2p-1$. 
\begin{align*}
\text{In}Q(v_1)=\langle x^p,x^{p-1}z,x^{p-2}z^3\dots,xz^{2p-3},z^{2p-1}\rangle\\
\text{In}Q(v_2)=\langle y^p,y^{p-1}z,y^{p-2}z^3\dots,yz^{2p-3},z^{2p-1}\rangle
\end{align*}
By lemma \ref{lemma7}, 
\begin{equation*}
\reg  H_0(J_\bullet)=\deg h_0(x,y)z^{\zeta_0}-3+(r+1)=2p+1-3+2p=4p-2=2r.
\end{equation*}
\end{itemize}
\end{proof}

\section*{Acknowledgement}
I would like to thank Michael DiPasquale for the discussion on Lemma \ref{lemma_fullrank}.
\bibliographystyle{plain}
\bibliography{Reference}
\end{document}